\documentclass[12pt]{amsart}

\usepackage[latin1]{inputenc}

\usepackage{amsfonts}
\usepackage{amssymb, amsmath}

\newtheorem{thm}{Theorem}
\newtheorem{lem}[thm]{Lemma}
\newtheorem{pro}[thm]{Proposition}

\title{Proximality and equidistribution on the Furstenberg boundary}
\author{A.~Gorodnik and F.~Maucourant}
\date{}
\thanks{The first author is partially supported by NSF grant 0400631.}
\address{Mathematics department\\ University of Michigan\\  Ann Arbor, MI 48109 USA}
\email{gorodnik@umich.edu}
\address{\'Ecole normale sup\'erieure de Lyon\\ Unit\'e de Math\'ematiques Pures et Appliqu\'ees\\
UMR CNRS 5669 \\ 46, all\'ee d'Italie\\  69364 Lyon Cedex 07 France}
\email{Francois.MAUCOURANT@umpa.ens-lyon.fr}

\begin{document}

\begin{abstract}
Let $G$ be a connected semisimple Lie group with finite center and without compact factors, $P$ a minimal parabolic subgroup of $G$,
and $\Gamma$ a lattice in $G$. We prove that every $\Gamma$-orbits in the Furstenberg boundary $G/P$
is equidistributed for the averages over Riemannian balls. The proof is based on the proximality of the action of $\Gamma$ on $G/P$.
\end{abstract}

\maketitle

\section{Introduction}
 
Let $G$ be a connected semisimple Lie group with finite center and without compact factor, and $\Gamma$ a lattice in $G$,
that is, a discrete subgroup of $G$ such that $\Gamma\backslash G$ has finite volume.
In this article we investigate the distribution of orbits of $\Gamma$ acting
on the Furstenberg boundary of $G$. Recall that the Furstenberg boundary can
be identified with the factor space $G/P$, where $P$ is a minimal parabolic subgroup of $G$.
It is known that every orbit of $\Gamma$ in $G/P$ is dense (see \cite{most}).
We show that orbits of $\Gamma$ are equidistributed with respect to the averages over Riemannian balls.

Since we study the action of a nonamenable group on a space without a finite invariant measure,
our result lies outside the scope of the classical ergodic theory.
The published results about distribution of dense orbits of nonamenable groups are limited to a few special examples.
Arnold and Krylov showed in \cite{ak} that dense orbits of groups generated by two rotations acting
on the $2$-dimensional sphere are equidistributed.
A similar problem was considered by Kazhdan in \cite{ka} where he studied the action of a group
generated by two affine isometries on the plane $\mathbb{R}^2$. 
Distribution of dense orbits of a lattice in $\hbox{\rm SL}(2,\mathbb{R})$
acting on $\mathbb{R}^2$ was investigated by Ledrappier \cite{l} and Nogueira \cite{no}.

Let $X$ be the symmetric space of $G$ equipped with a right invariant Riemannian metric $d$.
Note that $X$ can be identified with $L\backslash G$ for a maximal compact subgroup $L$ of $G$.

Fix $x,\tilde x\in X$ and denote by $K$ and $\tilde K$ the stabilizers of $x$ and $\tilde x$ respectively. Let $\nu$ and $\tilde\nu$ be the probability Haar measures on $K$ and $\tilde K$
and $m_{\tilde x}$ the harmonic measures at $\tilde x$ on $G/P$, that is, the unique $\tilde{K}$-invariant
probability measure on $G/P$.
For $S\subset G$ and  $T>0$, define 
\begin{eqnarray*}
S_T(\tilde x)&=&\{s\in S: d(x,\tilde xs)<T\},\\
S_T&=& S_T(x).
\end{eqnarray*}

Our main result is the following theorem.

\begin{thm}\label{th_main}
For every $f\in C(G/P)$, $\tilde x\in X$, and $y\in G/P$,
$$
\lim_{T\to\infty}\frac{1}{|\Gamma_T(\tilde x)|}\sum_{\gamma\in\Gamma_T(\tilde x)} f(\gamma y)= \int_{G/P} f dm_{\tilde x},
$$
Moreover, the convergence is uniform for $y\in G/P$.
\end{thm}

We remark that it was shown in \cite{em} (see also \cite{drs}) that 
\begin{equation}\label{eq_eskin}
|\Gamma_T(\tilde x)|\sim_{T\to\infty} \frac{\textrm{Vol}(G_T(\tilde x))}
{\textrm{Vol}(\Gamma\backslash G)},
\end{equation}
and the exact asymptotics of the volume $\textrm{Vol}(G_T(\tilde x))=\textrm{Vol}(G_T)$ as $T\to\infty$
was computed in \cite{kni}.

The first result in the direction of Theorem \ref{th_main} was established in \cite{mau},
where the case of the real hyperbolic spaces was considered. 
A different proof of Theorem \ref{th_main} is given in \cite{g_o}.
An advantage of the approach presented here is that it shows that the convergence is uniform.
While the proof in \cite{g_o} uses equidistribution of solvable flows on $\Gamma\backslash G$,
our proof is based on the strong proximality of the action of $G$ on $G/P$
(see Theorem \ref{th_proxy} below). This result is of independent interest, and
it might be useful for other applications.

Recall that an action of a group $H$ on a compact metric space $(Y,d)$ is called {\it proximal} if
for every $u,v\in Y$ there exists a sequence $\{h_n\}\subset H$ such that
$d(h_nu,h_nv)\to 0$ as $n\to\infty$. The fact that the action of $G$ on $G/P$ is proximal plays
important role in the study of random walks on $G$ (see, for example, \cite{furst}). It turns out that a typical sequence in $G$ acts on $G/P$ in proximal fashion.

\begin{thm}[\bf Strong proximality]\label{th_proxy}
Let $\mathcal{O}$ be neighborhood of the diagonal in $G/P\times G/P$ and $u,v\in G/P$. Then
$$
\lim_{T\to\infty}\frac{\hbox{\rm Vol}(\{g\in G_T(\tilde x): (gu,gv)\notin \mathcal{O}\})}{\hbox{\rm Vol}(G_T(\tilde x))}= 0
$$
and
$$
\lim_{T\to\infty}\frac{|\{\gamma\in \Gamma_T(\tilde x): (\gamma u,\gamma v)\notin \mathcal{O}\}|}{|\Gamma_T(\tilde x)|}= 0
$$
uniformly on $u,v$.
\end{thm}

In the case of the real hyperbolic space, Theorem \ref{th_proxy} was proved in \cite{mau}
using geometric methods.

\section{Proof of Theorem \ref{th_proxy}}

\subsection{Cartan decomposition}
Let $G=K_0\exp(\mathfrak{p})$ be the Cartan decomposition of $G$ and $A\subset \exp(\mathfrak{p})$
a split Cartan subgroup of $G$, that is, a maximal connected abelian subgroup in $\exp(\mathfrak{p})$.
We fix a system of positive roots $\Sigma^+$ on $\mathfrak{a}=\hbox{Lie}(A)$, and let 
$$
A^+=\{a\in A: \alpha(\log a)\ge 0\hbox{ for all }\alpha\in\Sigma^+\}
$$
denote the closed positive Weyl chamber in $A$. Then $G=KA^+K$, and
a Haar measure on $G$ can be given by
\begin{equation}\label{eq_cartan}
\int_{G}\psi(g)dg=\int_{K}\int_{A^+}\int_{K} \psi(k_1ak_2)\xi(\log a)d\nu(k_1)dad\nu(k_2),\quad \psi\in C_c(G),
\end{equation}
where $da$ denotes the Lebesgue measure on $A$,
$$
\xi(s)=\prod_{\alpha\in\Sigma^+} \sinh(\alpha(s))^{m_\alpha},\;\;\; s\in\mathfrak{a},
$$
and $m_\alpha$ denotes the dimension of the root space for the root $\alpha\in\Sigma^+$.

Let $\tilde g\in G$ be such that $x\tilde g=\tilde x$. Then $G=\tilde g^{-1} KA^+K$, $G_T(\tilde x)=\tilde g^{-1} KA^+_TK$, and
\begin{equation}\label{eq_cartan2}
\int_{G}\psi(g)dg=\int_{K}\int_{A^+}\int_{K} \psi(\tilde g^{-1}k_1ak_2)\xi(\log a)d\nu(k_1)dad\nu(k_2),\quad \psi\in C_c(G).
\end{equation}
In particular, it follows that
\begin{equation}\label{eq_gt}
\hbox{Vol}(G_T(\tilde x))=\hbox{Vol}(G_T)=\int_{A^+_T} \xi(\log a)da.
\end{equation}

\subsection{Reduction to maximal parabolics}
Fix a system of simple roots $$\Pi=\{\alpha_1,\ldots,\alpha_r\}\subset\Sigma^+.$$
Here $r=\dim A$ is the $\mathbb{R}$-rank of $G$. 
It is well-known that the closed subgroups of $G$ that contain $P$ are in
one-to-one correspondence with the subsets of $\Pi$ (see \cite[Sec.~1.2]{wa}). In particular, $P_i=P_{\{\alpha_i\}}$, $i=1,\ldots,r$,
are the maximal parabolic subgroups of $G$ and $$P=\bigcap_{i=1}^r P_i.$$
We consider the projection maps  
$$\pi_i : G/P\times G/P \rightarrow G/P_i\times G/P_i,\quad i=1,\ldots,r.$$
Let $\Delta$ and $\Delta_i$ denote the diagonals in $G/P\times G/P$ and $G/P_i\times G/P_i$ respectively.
Then $$\Delta=\bigcap_{i=1}^r \pi_i^{-1}(\Delta_i).$$ Since 
$$\prod_{i=1}^r\pi_i :\; G/P\times G/P \rightarrow \prod_{i=1}^r G/P_i\times G/P_i$$
is a continuous injective map from a compact space to a Hausdorff space, it is a homeomorphism onto its image.
It follows that for any neighborhood $\mathcal{O}$ of $\Delta$ in $G/P\times G/P$, there exist
neighborhoods $\mathcal{O}_i$ of $\Delta_i$ in $G/P_i\times G/P_i$ such that
$$
\mathcal{O}\supset \bigcap_{i=1}^r \pi_i^{-1}(\mathcal{O}_i).
$$
Then for every $(u,v)\in G/P\times G/P$,
$$
\{g\in G: g\cdot (u,v)\notin\mathcal{O}\}\subset \bigcup_{i=1}^r \{g\in G: g\cdot \pi_i(u,v)\notin\mathcal{O}_i\}.
$$
This inclusion shows that it suffices to prove Theorem \ref{th_proxy} under the assumption
that $P$ is a maximal parabolic subgroup of $G$. We keep this assumption until the end of this section.

\subsection{Dynamics on projective space}
By a result from \cite{Tits}, there is an irreducible representation $G\to\hbox{GL}(V)$ such that
the highest weight space is one-dimensional, and the stabilizer of this space is $P$.
We consider the induced action of $G$ on the projective space $\mathbb{P}(V)$, and let $w^+\in \mathbb{P}(V)$
be the direction of the highest weight space. The map $g\mapsto gw^+$ defines an embedding of $G/P$ in $\mathbb{P}(V)$.
Note that if $\lambda$ is the highest weight, the other weights of the representation are of the form
$\lambda-\sum_{\alpha\in\Sigma^+} n_\alpha \alpha$ for integers $n_\alpha\ge 0$.
We denote by $V^<$ the sum of all root spaces with weights other than $\lambda$.     
We fix a $K$-invariant scalar product on $V$,  which gives rise to a metric $d$ on $\mathbb{P}(V)$, which is $K$-invariant. Put $\tilde d(w_1,w_2)= d(\tilde gw_1,\tilde gw_2)$. 
Let $V^<_\varepsilon$ be the open $\varepsilon$-neighborhood of $V^<$ in $\mathbb{P}(V)$ with respect to the
metric $\tilde d$.

For $w \in  \mathbb{P}(V)$ and $\tau>0$, define
$$
K_{\tau}(w)=\{ k \in K : kw \notin V^<_\tau\}.
$$
   \begin{lem}\label{l_k_e}
   For every $w\in G\cdot w^+$,
   $$\lim_{\tau \rightarrow 0^+} \nu(K-K_{\tau}(w)) =0.$$
   \end{lem}
   \begin{proof}
   It follows from the Iwasawa decomposition that $G\cdot w^+=K\cdot w^+$.
   Thus, without loss  of generality, we may assume that $w=w^+$.
   By the continuity of the measure, it suffices to prove that
   $$
   \nu(\{k\in K: kw^+\in V^<\})=0.
   $$
   Suppose that this is false. For a subspace $W$ of $V$, define
   $$
   K_W=\{k\in K: kw^+\in W\}.
   $$ 
   Let $W$ be a minimal subspace of $V^<$ such that $\nu(K_W)>0$. We claim that $\hbox{Stab}_{K}(W)=K$.
   If $\hbox{Stab}_{K}(W)$ has infinite index in $K$, then there exist $k_i\in K$, $i\ge 1$, such that
   $k_iW\ne k_j W$ for $i\ne j$. Since all sets $k_iK_W\subset K$, $i\ge 1$, have the same positive measure,
   it follows that for some $i\ne j$, $k_iK_W\cap k_jK_W$ has positive measure.
   Then $k_j^{-1}k_iK_W\cap K_W$ has positive measure too, and for $k\in k_j^{-1}k_iK_W\cap K_W$, $$kw^+\in k_j^{-1}k_iW\cap W.$$
   Since $k_j^{-1}k_iW\cap W$ is a proper subspace of $W$, this contradicts the choice of $W$. Thus, 
   $\hbox{Stab}_{K}(W)$ is a closed subgroup of finite index in $K$. Since $K$ is connected,
   it follows that $K=\hbox{Stab}_{K}(W)$. Then $w^+\in K_W^{-1}W\subset V^<$. This contradiction proves the lemma. 
    \end{proof}

 We consider the sets
 \begin{eqnarray}
 A_T^\eta&=&\{a\in A_T: \alpha(\log a)\ge \eta\hbox{ for all }\alpha\in\Sigma^+\},\nonumber\\
 G_{T,\varepsilon}(u,v)&=&\{g \in G_T(\tilde x) : \tilde d(gu,gv)>\varepsilon\},\label{eq_ate}\\
 \Omega^\eta_{T,\tau}(u,v)&=&\tilde g^{-1}K A^\eta_{T} (K_{\tau}(u)\cap K_{\tau}(v))\nonumber
 \end{eqnarray}
 defined for $T,\eta,\tau,\varepsilon>0$ and $u,v\in\mathbb{P}(V)$.
   
   \begin{lem}\label{l_last}
   For every $\varepsilon>0$ and $\tau>0$, there exists $\eta>0$ such that for every $T>0$
     and $u,v \in G\cdot w^+$,
     \begin{equation}\label{eq_empty}
     \Omega^\eta_{T,\tau}(u,v) \cap G_{T,\varepsilon}(u,v) = \emptyset.
     \end{equation}
   \end{lem}
   \begin{proof}
     Note that an element $a \in A^\eta_{T}$ acts by diagonal transformations
     on $V$ with respect to some fixed basis, and the eigenvalue associated to the vector $w^+$
     is at least $e^\eta$ times greater than the other eigenvalues. Therefore, for all $w\notin V_\tau^<$ and
     sufficiently large $\eta$ (depending only on $\tau$ and $\varepsilon$), we have $d(aw,w^+)<\varepsilon/2$ when $a \in A^\eta_{T}$.
     Thus, for
     $$\tilde g^{-1}k_1 a k_2 \in \Omega^\eta_{T, \tau}(u,v)=\tilde g^{-1}K A^\eta_{T} (K_{\tau}(u)\cap K_{\tau}(v)),$$
     we have
     $$\tilde d(\tilde g^{-1} k_1a k_2 u,\tilde g^{-1} k_1a k_2 v)=d(a k_2 u, a k_2v)\leq d(a k_2 u ,w^+)+d(a k_2 v ,w^+)<\varepsilon,$$
     This proves the lemma.
   \end{proof}

\subsection{Completion of the proof}   
By (\ref{eq_cartan2}),
\begin{equation}\label{eq_Ht}
\hbox{Vol}(\Omega^\eta_{T,\tau}(u,v))=\left(\int_{A^\eta_{T}}\xi (\log a) da\right)\cdot \nu(K_{\tau}(u)\cap K_{\tau}(v)).
\end{equation}
Let $\varepsilon,\delta\in (0,1)$. Using Lemma \ref{l_k_e}, we choose $\tau>0$ such that 
$$
\nu(K_{\tau}(u)\cap K_{\tau}(v))> 1-\delta.
$$
Let $\eta>0$ be as Lemma \ref{l_last}. By Lemma \ref{l_a_eta}(a), for sufficiently large $T$,
$$
\int_{A^\eta_{T}}\xi (a) da\ge (1-\delta) \int_{A^+_{T}}\xi (\log a) da.
$$
Thus, it follows from (\ref{eq_gt}) and (\ref{eq_Ht}) that
$$
\hbox{Vol}(\Omega^\eta_{T,\tau}(u,v))\ge (1-\delta)^2\hbox{Vol}(G_T(\tilde x)).
$$
for sufficiently large $T>0$. Therefore, by (\ref{eq_empty}),
$$
\hbox{Vol}(G_{T,\varepsilon}(u,v))\le (1-(1-\delta)^2)\hbox{Vol}(G_T(\tilde x))
$$
for all $\delta \in (0,1)$ and sufficiently large $T>0$.
Since the sets $$\{(g_1P,g_2P): \tilde d(g_1w^+,g_2w^+)<\varepsilon\},\quad \varepsilon >0,$$ form a base of
the neighborhoods of the diagonal in $G/P\times G/P$,
this proves the first part of Theorem \ref{th_proxy}.

To prove the second part of Theorem \ref{th_proxy}, we choose a neighborhood $\mathcal{P}$
of $e$ in $G$ and a neighborhood $\mathcal{Q}$ of the diagonal in $G/P\times G/P$ such that
\begin{eqnarray}
\mathcal{P}^{-1}\mathcal{P}\cap \Gamma&=&\{e\},\label{eq_1}\\
\mathcal{P}^{-1}\cdot\mathcal{Q}&\subset& \mathcal{O},\label{eq_2}\\
\mathcal{P}\cdot G_T(\tilde x)&\subset& G_{T+c}(\tilde x).\label{eq_3}
\end{eqnarray}
for fixed $c>0$ and all $T>0$.
Here we use that $\Gamma$ is discrete, the space $G/P$ is compact, and the metric
on the symmetric space is uniformly continuous.
By (\ref{eq_2}), for every $\gamma\in\Gamma$ such that $\gamma\cdot (u, v)\notin \mathcal{O}$, 
we have $\mathcal{P}\gamma\cdot (u,v)\cap \mathcal{Q}=\emptyset$. 
Thus, using (\ref{eq_3}), we deduce that
$$
\mathcal{P}\cdot \{\gamma\in\Gamma_T(\tilde x):\gamma\cdot (u, v)\notin \mathcal{O}\}
\subset \{g\in G_{T+c}(\tilde x): g\cdot (u, v)\notin \mathcal{Q}\}.
$$
Then by (\ref{eq_1}), $\mathcal{P}\gamma_1\cap \mathcal{P}\gamma_2=\emptyset$ for $\gamma_1,\gamma_2\in\Gamma$, $\gamma_1\ne \gamma_2$, and
\begin{eqnarray*}
|\{\gamma\in\Gamma_T(\tilde x):\gamma\cdot (u, v)\notin \mathcal{O}\}|&\le&
\frac{1}{\hbox{Vol}(\mathcal{P})}\hbox{Vol}(\{g\in G_{T+c}(\tilde x): g\cdot (u, v)\notin \mathcal{Q}\})\\
&=& o(\hbox{Vol}(G_{T+c}(\tilde x)))
\end{eqnarray*}
as $T\to\infty$. Now the second statement of Theorem \ref{th_proxy} follows from Lemma \ref{l_a_eta}(d) and (\ref{eq_eskin}).

\section{Equidistribution on $\Gamma\backslash G$}

Recall that $K$ is a maximal compact subgroups of $G$, and $\nu$
is the probability Haar measure on $K$.
Denote by $\varrho$ a right Haar measure on the minimal parabolic subgroup $P$. For a suitable normalization of $\varrho$,
the Haar measure on $G$ is given by
\begin{equation}\label{eq_iwasawa}
\int_G \psi(g) dg=\int_{K}\int_P \psi(kp)d\varrho(p)d\nu(k),\quad \psi\in C_c(G).
\end{equation}
We also define a measure $\mu$ on $G$ by
\begin{equation}\label{eq_mu}
\int_G \psi(g) d\mu(g)=\int_{K}\int_P \psi(kp^{-1})d\varrho(p)d\nu(k),\quad \psi\in C_c(G).
\end{equation}
Note that $\mu$ is left $K$-invariant.

The first step in the proof of Theorem \ref{th_main} is the following result.
\begin{pro} \label{p_gg}
For every $\Psi\in C_c(\Gamma\backslash G)$ and $z\in \Gamma\backslash G$,
$$
\lim_{T\to\infty}\frac{1}{\mu(G_T)}\int_{G_T} \Psi(zg)d\mu(g)=\frac{1}{\hbox{\rm Vol}(\Gamma\backslash G)}\int_{\Gamma\backslash G} \Psi\;dg
$$
where $G_T=\{g\in G: d(x,xg)<T$\}.
\end{pro}
Proposition \ref{p_gg} is a consequence of the equidistribution of translates of $K$ in
$\Gamma\backslash G$ proved by Eskin and McMullen in \cite{em} (see also \cite{sh} for a more general result).
They showed that for every strongly divergent sequence $\{g_n\}\subset G$,
\begin{equation}\label{eq_K}
\lim_{n\to\infty} \int_{K} \Psi(zkg_n)d\nu(k)=\frac{1}{\hbox{\rm Vol}(\Gamma\backslash G)}\int_{\Gamma\backslash G} \Psi\;dg.
\end{equation}
Recall that a sequence $\{g_n\}\subset G$ is {\it strongly divergent} if the projection of $\{g_n\}$
on every noncompact simple factor of $G$ is divergent.
Note that (\ref{eq_K}) was proved in \cite{em} under the condition that the lattice $\Gamma$ is irreducible.
Since the proof of (\ref{eq_K}) is based on mixing properties of the action of $G$ on $\Gamma\backslash G$,
it is applicable to the case of a reducible lattice $\Gamma$ provided that the sequence $\{g_n\}$  is strongly divergent.

Denote by $\pi_i:G\to G_i$, $i=1,\ldots,s$, the projections of $G$ onto its simple factors.
Let $C_{i,j}\subset G_i$, $j\ge 1$, be an increasing sequence of compact subsets such that $G_i=\cup_{j\ge 1} C_{i,j}$.
Define
\begin{equation}\label{eq_gtn}
G_{T,n}=G_T- \bigcup_{1\le i \le s} \pi_i^{-1}(C_{i,n}).
\end{equation}

\begin{lem}\label{l_gtn}
For every $n\ge 1$, $\mu(G_{T,n})\sim \mu(G_T)$ as $T\to\infty$.
\end{lem}

\begin{proof}
It suffices to show that for every $i=1,\ldots,s$ and $n\ge 1$,
$$
\mu(G_T\cap \pi_i^{-1}(C_{i,n}))=o(\mu(G_T))\quad\hbox{as}\quad T\to\infty.
$$
Fix $i=1,\ldots,s$ and $n\ge 1$.
Note that $G=DH$, where $D$ and $H=\hbox{ker}(\pi_i)$ are normal connected
semisimple Lie subgroups with finite centers, and $D$ and $H$ commute.
We have $\pi_i^{-1}(C_{i,n})=D_{i,n}H$ for some compact set $D_{i,n}\subset D$.
There is a constant $\delta>0$ such that 
\begin{equation}\label{eq_DH}
D_{i,n}H_{T-\delta}\subset (D_{i,n}H)_T\subset D_{i,n}H_{T+\delta}\;\;\hbox{ for all $T>0$.}
\end{equation}
We define measures $\mu_D$ and $\mu_H$ for the groups $D$ and $H$ respectively as in (\ref{eq_mu}).
With appropriate normalization, $\mu=\mu_D\otimes\mu_H$. Thus, it follows from (\ref{eq_DH}) that
\begin{equation}\label{eq_Ha}
\mu(G_T\cap \pi_i^{-1}(C_{i,n}))=\mu((D_{i,n}H)_T)\ll \mu_H(H_{T+\delta}).
\end{equation}
Since $G_T=KP_T$ and $P_T^{-1}=P_T$, using (\ref{eq_iwasawa}) and (\ref{eq_mu}), we conclude that
\begin{equation}\label{eq_mu_vol}
\mu(G_T)=\varrho(P_T^{-1})=\varrho(P_T)=\hbox{Vol}(G_T).
\end{equation}
Similarly, $H=LQ_T$ where $L$ is a maximal compact subgroup of $H$ contained in $K$, and $Q$ is a minimal
parabolic subgroup of $H$. As in (\ref{eq_mu_vol}), we deduce that $\mu_H(H_T)=\hbox{Vol}_H(H_T)$.
By (\ref{eq_Ha}), it is sufficient to show that
\begin{equation}\label{eq_HTa}
\hbox{Vol}_H(H_{T+\delta})=o(\hbox{Vol}(G_{T}))\quad\hbox{as}\quad T\to\infty.
\end{equation}
Note that with appropriate normalization the Haar measure on $G$ is the product of Haar measures on $D$ and $H$.
Without loss of generality, $\hbox{Vol}_D(D_{i,n})>0$. Then by (\ref{eq_DH}),
$$
\hbox{Vol}_H(H_{T+\delta})\ll \hbox{Vol}(D_{i,n}H_{T+\delta})\le \hbox{Vol}((D_{i,n}H)_{T+2\delta}).
$$
Let $G_T^\eta$ be defined as in (\ref{eq_last}). Since the set $D_{i,n}$ is compact, there
exists $\eta>0$ such that
$$
(D_{i,n}H)_{T+2\delta}\subset G_{T+2\delta}- G^\eta_{T+2\delta}.
$$
Thus, (\ref{eq_HTa}) follows from Lemma \ref{l_a_eta}(b).
\end{proof}

\begin{proof}[Proof of Proposition \ref{p_gg}]
The map $K\times A^+\times K\to G$ is a diffeomorphism on an open set of full measure.
Since the measure $\mu$ is left $K$-invariant and smooth, for some $\sigma\in C(A^+\times K)$,
$$
\int_G \psi (g)d\mu(g)=\int_{K}\int_{A^+}\int_{K} \psi (k_1ak_2)\sigma(a,k_2) d\nu(k_1)dad\nu(k_2),\quad \psi\in C_c(G).
$$

Let $G_{T,n}$ be defined as in (\ref{eq_gtn}), and it is $K$-bi-invariant (equivalently,
all $C_{i,j}$ are $\pi_i(K)$-bi-invariant). Then 
$$
G_{T,n}=KA^+_{T,n}K\;\;\hbox{and}\;\; 
\mu(G_{T,n})=\int_{K}\int_{A^+_{T,n}}\sigma(a,k_2)dad\nu(k_2),
$$
where $A^+_{T,n}=G_{T,n}\cap A^+$.

Let $\varepsilon>0$. By (\ref{eq_K}), 
$$
\left|\int_{K} \Psi(zk_1ak_2)d\nu(k_1)-\frac{1}{\hbox{\rm Vol}(\Gamma\backslash G)}\int_{\Gamma\backslash G} \Psi\;dg \right|<\varepsilon
$$
for $a\in A^+_{T,n}$ and $k_2\in K$ when $n>n_0(\varepsilon)$. Thus, for $n>n_0(\varepsilon)$,
\begin{eqnarray}\label{eq_g_gt}
&&\left|\int_{G_{T,n}} \Psi(zg)d\mu(g)-\frac{\mu(G_{T,n})}{\hbox{\rm Vol}(\Gamma\backslash G)}\int_{\Gamma\backslash G} \Psi\;dg \right|\\
&=&\left|\int_{K}\int_{A^+_{T,n}}\int_{K} \Psi(zk_1ak_2)d\nu(k_1)\sigma(a,k_2)dad\nu(k_2)\right.\nonumber \\
&-&\left.\frac{\mu(G_{T,n})}{\hbox{\rm Vol}(\Gamma\backslash G)}\int_{\Gamma\backslash G} \Psi\;dg \right|
\le \int_{K}\int_{A^+_{T,n}}\left|\int_{K}\Psi(zk_1ak_2)d\nu(k_1)\right.\nonumber \\
&-& \left. \frac{1}{\hbox{\rm Vol}(\Gamma\backslash G)}\int_{\Gamma\backslash G} \Psi\;dg\right|\sigma(a,k_2)dad\nu(k_2)<\varepsilon \mu(G_{T,n}).\nonumber
\end{eqnarray}
By Lemma \ref{l_gtn}, for every $n\ge 1$,
$$
\int_{G_{T}} \Psi(zg)d\mu(g)=\int_{G_{T,n}} \Psi(zg)d\mu(g)+o(\mu(G_{T,n}))
$$
as $T\to\infty$. Thus, it follows from (\ref{eq_g_gt}) that
$$
\limsup_{T\to\infty}\left|\frac{1}{\mu(G_T)}\int_{G_{T}} \Psi(zg)d\mu(g)-\frac{1}{\hbox{\rm Vol}(\Gamma\backslash G)}\int_{\Gamma\backslash G} \Psi\;dg \right|<\varepsilon
$$
for every $\varepsilon>0$. This proves the proposition.
\end{proof}

\section{Equidistribution on average}

In this section we prove that Theorem \ref{th_main} holds ``on average''.
In the case of hyperbolic spaces, the following proposition is a consequence
 of the work
 of Roblin \cite{ro}.
\begin{pro}\label{p_average}
For every $f\in C(G/P)$ and $y\in G/P$,
$$
\lim_{T\to\infty}\frac{1}{|\Gamma_T(\tilde x)|}\sum_{\gamma\in\Gamma_T(\tilde x)} \int_{K} f(\gamma k y)d\nu(k)=\int_{G/P} f dm_{\tilde x}
$$
where $\Gamma_T(\tilde x)=\{\gamma\in \Gamma: d(x,\tilde x\gamma)<T\}$.
\end{pro}
\begin{proof}
There exists $\tilde p\in P$ such that $\tilde x=x\tilde p$. Then $\tilde K=\tilde p^{-1}K\tilde p$, and
it follows from (\ref{eq_iwasawa}) that
\begin{equation}\label{eq_iwasawa2}
\int_G \psi(g) dg=\int_{\tilde K}\int_P \psi(k\tilde p^{-1}p)d\varrho(p)d\tilde \nu(k),\quad \psi\in C_c(G).
\end{equation}

Without loss of generality, $f\ge 0$, and since $G=KP$, we may assume that $y=eP$.

Let $\varepsilon>0$, $\mathcal{O}_\varepsilon=\{z\in X: d(x,z)<\varepsilon\}$, and $\phi_\varepsilon\in C_c(X)$ such that 
$$
\phi_\varepsilon\ge 0,\quad \hbox{supp}(\phi_\varepsilon)\subset \mathcal{O}_\varepsilon, \quad \int_P \phi_\varepsilon(xp^{-1}) d\varrho(p)=1.
$$
Since $X=\tilde xP$ and $\varrho$ is right invariant, it follows that 
\begin{equation}\label{eq_int_P}
\int_P \phi_\varepsilon(zp^{-1}) d\varrho(p)=1\quad\hbox{ for every $z\in X$.}
\end{equation}
Let $$\psi_\varepsilon(g)=f(gP)\phi_\varepsilon(\tilde xg),\quad g\in G.$$
Clearly, $\psi_\varepsilon\in C_c(G)$ and 
$$
\Psi_\varepsilon(\Gamma g)\stackrel{def}{=}\sum_{\gamma\in\Gamma} \psi_\varepsilon(\gamma g)\in C_c(\Gamma\backslash G).
$$
By Proposition \ref{p_gg},
\begin{equation}\label{eq_mu_lim0}
\lim_{T\to\infty}\frac{1}{\mu(G_T)}\sum_{\gamma\in \Gamma} \int_{G_T} \psi_\varepsilon(\gamma g)d\mu(g)=\frac{1}{\hbox{Vol}(\Gamma\backslash G)}\int_{\Gamma\backslash G} \Psi_\varepsilon(\Gamma g) dg
\end{equation}
and by (\ref{eq_iwasawa2}),
\begin{eqnarray*}
\hbox{Vol}(\Gamma\backslash G)\int_{\Gamma\backslash G} \Psi_\varepsilon(\Gamma g) dg &=& \int_G \psi_\varepsilon(g)dg 
=\int_{\tilde K} f(kP)d\tilde\nu(k)\cdot\int_P \phi_\varepsilon(\tilde x\tilde p^{-1}p)d\varrho(p)\\
&=&\int_{G/P} f dm_{\tilde x}\cdot\int_P \phi_\varepsilon(xp)d\varrho(p).\nonumber 
\end{eqnarray*}
Denote by $\delta$ the modular function of $P$. By (\ref{eq_int_P}),
\begin{eqnarray*}
\left| \int_P \phi_\varepsilon(xp)d\varrho(p)-1\right|&=&\left| \int_P \phi_\varepsilon(xp^{-1})(\delta(p)-1)d\varrho(p)\right|\\
&\le&\max\{|\delta(p)-1|: xp^{-1}\in\mathcal{O}_\varepsilon\}.
\end{eqnarray*}
The sets $\{p\in P: xp^{-1}\in\mathcal{O}_\varepsilon\}$, $\varepsilon>0$, form a base of neighborhoods of $P\cap K$ in $P$.
Since $\delta|_{P\cap K}=1$ and $P\cap K$ is compact, 
$$
\max\{|\delta(p)-1|: xp^{-1}\in\mathcal{O}_\varepsilon\}\to 0\quad\hbox{ as $\varepsilon\to 0^+$.}
$$
Thus, it follows from (\ref{eq_mu_lim0}) that
\begin{equation}\label{eq_mu_lim}
\lim_{\varepsilon\to 0^+}\lim_{T\to\infty}\frac{1}{\mu(G_T)}\sum_{\gamma\in \Gamma} \int_{G_T} \psi_\varepsilon(\gamma g)d\mu(g)=\int_{G/P} f dm_{\tilde x}.
\end{equation}
Since $G_T=KP_T$,
\begin{eqnarray*}
&&\sum_{\gamma\in \Gamma} \int_{G_T} \psi_\varepsilon(\gamma g)d\mu(g)\\
&\stackrel{(\ref{eq_mu})}{=}&\sum_{\gamma\in \Gamma}\int_{K\times P_T^{-1}} \psi_\varepsilon(\gamma kp^{-1})d\nu(k)d\varrho(p)\\
&=&\sum_{\gamma\in \Gamma} \int_{K} f(\gamma kP)
\left(\int_{P_T^{-1}} \phi_\varepsilon(\tilde x\gamma kp^{-1})d\varrho(p)\right)d\nu(k).
\end{eqnarray*}
For $\gamma\in \Gamma-\Gamma_{T+\varepsilon}(\tilde x)$, $k\in K$, and $p\in P_T^{-1}$,
$$
d(x,\tilde x\gamma kp^{-1})=d(xpk^{-1},\tilde x\gamma)\ge d(x,\tilde x\gamma)-d(x, xpk^{-1})\ge \varepsilon.
$$
This implies that $\int_{P_T^{-1}} \phi_\varepsilon(\tilde x\gamma kp^{-1})d\varrho(p)=0$ for $\gamma\in \Gamma-\Gamma_{T+\varepsilon}(\tilde x)$. Thus,
\begin{eqnarray*}
&&\sum_{\gamma\in \Gamma} \int_{G_T} \psi_\varepsilon(\gamma g)d\mu(g)\\
&=&\sum_{\gamma\in \Gamma_{T+\varepsilon}(\tilde x)} \int_{K} f(\gamma kP)
\left(\int_{P_T^{-1}} \phi_\varepsilon(\tilde x\gamma kp^{-1})d\varrho(p)\right)d\nu(k)\\
&\le&
\sum_{\gamma\in \Gamma_{T+\varepsilon}(\tilde x)} \int_{K} f(\gamma kP)
\left(\int_{P} \phi_\varepsilon(\tilde x\gamma kp^{-1})d\varrho(p)\right)d\nu(k)\\
&\stackrel{(\ref{eq_int_P})}{=}&\sum_{\gamma\in \Gamma_{T+\varepsilon}(\tilde x)} \int_{K} f(\gamma kP)d\nu(k).
\end{eqnarray*}
Combining (\ref{eq_mu_lim}), (\ref{eq_mu_vol}), (\ref{eq_eskin}) and Lemma \ref{l_a_eta}(c), we deduce that
$$
\liminf_{T\to\infty}\frac{1}{|\Gamma_T(\tilde x)|}\sum_{\gamma\in \Gamma_{T}(\tilde x)} \int_{K} f(\gamma kP)d\nu(k)\ge \int_{G/P} f dm_{\tilde x}.
$$

On the other hand, for $\gamma\in \Gamma_{T-\varepsilon}(\tilde x)$, $k\in K$, and $p\in P$ such that
$d(x,\tilde x\gamma kp^{-1})<\varepsilon$,
$$
d(x,xp^{-1})\le d(x,\tilde x\gamma kp^{-1})+ d(xp^{-1},\tilde x\gamma k p^{-1})<T.
$$
This shows that for $\gamma\in \Gamma_{T-\varepsilon}(\tilde x)$,
$$
\int_{P^{-1}_T} \phi_\varepsilon(\tilde x\gamma kp^{-1})d\varrho(p)=\int_{P} \phi_\varepsilon(\tilde x\gamma kp^{-1})d\varrho(p)\stackrel{(\ref{eq_int_P})}{=}1.
$$
Hence,
\begin{eqnarray*}
&&\sum_{\gamma\in \Gamma} \int_{G_T} \psi_\varepsilon(\gamma g)d\mu(g)\\
&\ge&
\sum_{\gamma\in \Gamma_{T-\varepsilon}(\tilde x)} \int_{K} f(\gamma kP)
\left(\int_{P_T^{-1}} \phi_\varepsilon(\tilde x\gamma kp^{-1})d\varrho(p)\right)d\nu(k)\\
&=&\sum_{\gamma\in \Gamma_{T-\varepsilon}(\tilde x)} \int_{K} f(\gamma kP)d\nu(k).
\end{eqnarray*}
By (\ref{eq_mu_lim}), (\ref{eq_mu_vol}), (\ref{eq_eskin}), and Lemma \ref{l_a_eta}(c),
$$
\limsup_{T\to\infty}\frac{1}{|\Gamma_T(\tilde x)|}\sum_{\gamma\in \Gamma_{T}(\tilde x)} \int_{K} f(\gamma kP)d\nu(k)\le \int_{G/P} f dm_{\tilde x}.
$$
This proves the proposition.
\end{proof}

\section{Proof of Theorem \ref{th_main}}

Now the proof can be completed using the argument from \cite{mau}.
Let $\varepsilon>0$. Since the space $G/P\times G/P$ is compact, there exists a neighborhood $\mathcal{O}$
of the diagonal in $G/P\times G/P$ such that for every $(z_1,z_2)\in \mathcal{O}$, we have
$|f(z_1)-f(z_2)|<\varepsilon$. Then for every $k\in K$,
\begin{eqnarray*}
&&\left|\sum_{\gamma\in\Gamma_T(\tilde x)} f(\gamma y)-\sum_{\gamma\in\Gamma_T(\tilde x)} f(\gamma k y)\right|\\
&\le& \sum_{\gamma\in\Gamma_T(\tilde x): (\gamma y,\gamma ky)\in\mathcal{O}} |f(\gamma y)-f(\gamma ky)|
+ \sum_{\gamma\in\Gamma_T(\tilde x): (\gamma y,\gamma ky)\notin\mathcal{O}} |f(\gamma y)-f(\gamma ky)|\\
&\le& \varepsilon |\Gamma_T(\tilde x)|+ 2 \sup |f|\cdot 
|\{\gamma\in\Gamma_T(\tilde x): (\gamma y,\gamma ky)\notin\mathcal{O}\}|.
\end{eqnarray*}
Thus, it follows from Theorem \ref{th_proxy} that
$$
\lim_{T\to\infty} \frac{1}{|\Gamma_T(\tilde x)|}\left|\sum_{\gamma\in\Gamma_T(\tilde x)} f(\gamma y)-\sum_{\gamma\in\Gamma_T(\tilde x)} f(\gamma k y)\right|=0
$$
for all $k\in K$. Hence, by the dominated convergence theorem,
$$
\lim_{T\to\infty} \left|\frac{1}{|\Gamma_T(\tilde x)|}\sum_{\gamma\in\Gamma_T(\tilde x)} f(\gamma y)-\frac{1}{|\Gamma_T(\tilde x)|}\sum_{\gamma\in\Gamma_T(\tilde x)}\int_{K} f(\gamma k y)d\nu(k)\right|=0.
$$
Finally, Theorem \ref{th_main} follows from Proposition \ref{p_average}.

\section{Appendix: volume estimates}

In this section, we give proofs of volume estimates, which are used in Theorems \ref{th_main} and \ref{th_proxy}.
There are other ways to establish these volume estimates. For example, one can use the exact asymptotics of
the volume of Riemannian balls from \cite{kni} (see also \cite{g_o}).
We present a straightforward proof that does not use asymptotics.

Let $\mathfrak{a}$ be the Lie algebra of the Cartan subgroup $A$ and $\mathfrak{a}^+$
the positive Weyl chamber with respect to the root system $\Sigma^+$.
The Riemannian metric defines a scalar product on $\mathfrak{a}$ and, by duality, on the dual space of $\mathfrak{a}$.
For $\alpha\in\Sigma^+$, we denote by $m_\alpha$
the dimension of the corresponding root space and put $\rho=\frac{1}{2}\sum_{\beta \in \Sigma^+} m_\beta \beta$.

   \begin{lem}\label{l_max}
   The maximum of $\rho$ on $\{a\in\mathfrak{a}:\|a\|\le 1\}$
   is achieved at a unique point in the interior of $\mathfrak{a}^+$.
   \end{lem}
   \begin{proof}
   Since the set $\{a\in\mathfrak{a}:\|a\|=1\}$ is strictly convex, it is clear that
   the point of maximum is unique.
   It is sufficient to show that $(\rho,\alpha)>0$ for every $\alpha\in\Sigma^+$.
   Denote by $\sigma_\alpha$ the reflection with respect to the hyperplane $\{\alpha=0\}$.
   The map $\sigma_\alpha$ permutes the elements of the set $\Sigma^+-\{\alpha,2\alpha\}$ and
   $\sigma_\alpha(\alpha)=-\alpha$. Since $m_{\sigma_\alpha(\beta)}=m_\beta$, we have
    $$\sigma_\alpha(\rho)=\rho-2m_\alpha \alpha-4m_{2\alpha}\alpha.$$
   Thus, 
$$(\rho,\alpha)=(\sigma_\alpha(\rho),\sigma_\alpha(\alpha))=2m_\alpha(\alpha,\alpha)+4m_{2\alpha}(\alpha,\alpha)-(\rho,\alpha)$$
    and $(\rho,\alpha)=(m_\alpha+2m_{2\alpha})(\alpha,\alpha)$ is positive.
    \end{proof}

For $T,\eta>0$, define
\begin{eqnarray}
A_{T}^{\eta}&=&\{a \in A_T \, : \alpha(\log a)\ge \eta\hbox{ for all } \alpha\in\Sigma^+ \}\nonumber\\
&=& \{a \in A \, : \|\log a\|<T,\,\alpha(\log a)\ge \eta\hbox{ for all } \alpha\in\Sigma^+ \},\label{eq_last}\\
G_T^{\eta}&=&K A_{T}^{\eta}K.\nonumber
\end{eqnarray}

  \begin{lem}\label{l_a_eta} For every $\eta>0$,
\begin{enumerate}  
\item[(a)]
    $$\int_{A_{T}^{\eta}} \xi(\log a)da \sim_{T\to\infty} \int_{A^+_T} \xi(\log a)da,$$
\item[(b)] 
$$\hbox{\rm Vol}(G^\eta_T)\sim_{T\to\infty} \hbox{\rm Vol}(G_T),$$
\item[(c)]
$$
\liminf_{\varepsilon\to 0^+}\left( \limsup_{T\to\infty} \frac{\hbox{\rm Vol}(G_{T+\varepsilon})}{\hbox{\rm Vol}(G_{T})}\right)=1,
$$
\item[(d)]
$$
\hbox{\rm Vol}(G_{T+\eta})\ll \hbox{\rm Vol}(G_T).
$$
\end{enumerate}
   \end{lem}
   \begin{proof}
We have
\begin{equation}\label{eq_sum}
\int_{\mathfrak{a}^+_T} \xi(a)da=2^{-|\Sigma^+|}\sum_{i\in I} \int_{\mathfrak{a}_T^+} e^{\lambda_i(a)}da
\end{equation}
where $\lambda_i$'s the characters of the form $2\rho-\sum_{\alpha\in\Sigma^+} n_\alpha\alpha$ for some $n_\alpha\ge 0$.
Let
\begin{eqnarray*}
\delta&=&\max\{2\rho(a): a\in\mathfrak{a}^+_1\},\\
\delta_i&=&\max\{\lambda_i(a): a\in\mathfrak{a}^+_1\},\;\;\; i\in I,\\
\delta_\alpha&=&\max\{2\rho(a): a\in\mathfrak{a}^+_1,\alpha(a)=0\},\;\;\;\alpha\in\Sigma^+.
\end{eqnarray*}
It follows from Lemma \ref{l_max} that for $\lambda_i\ne 2\rho$ and $\alpha\in\Sigma^+$,
$\delta>\max\{\delta_i,\delta_\alpha\}$. 
Thus,
\begin{equation}\label{eq_lambda}
\int_{\mathfrak{a}^+_T} e^{\lambda_i(a)}da\le \hbox{Vol}(\mathfrak{a}_T^+) e^{\delta_i T}\ll T^{r} e^{\delta_i T}
\end{equation}
where $r=\dim \mathfrak{a}$.
Let $\varepsilon>0$ be such that 
$$
\delta-\varepsilon>\max\{\delta_i,\delta_\alpha: i\in I,\alpha\in\Sigma^+\}.
$$
Then
\begin{eqnarray}
\int_{\mathfrak{a}^+_T} e^{2\rho(a)}da &=& T^r\int_{\mathfrak{a}_1^+} e^{2T\rho(a)}da\label{eq_up}\\
&\ge& T^r e^{(\delta-\varepsilon)T}\hbox{Vol}(\{a\in \mathfrak{a}_1^+: 2\rho(a)\ge \delta-\varepsilon\})\gg T^r e^{(\delta-\varepsilon)T}.\nonumber
\end{eqnarray}
Combining (\ref{eq_sum}), (\ref{eq_lambda}), and (\ref{eq_up}), we deduce that
\begin{equation}\label{eq_low}
\int_{\mathfrak{a}^+_T} \xi(a)da\gg T^r e^{(\delta-\varepsilon)T}.
\end{equation}
On the other hand, for $\alpha\in\Sigma^+$,
\begin{eqnarray*}
\int_{\mathfrak{a}^+_T\cap\{\alpha<\eta\}} \xi(a)da&\le& \int_{\mathfrak{a}^+_T\cap\{\alpha<\eta\}} e^{2\rho(a)}da
\ll \int_{\mathfrak{a}^+_T\cap\{\alpha=0\}} e^{2\rho(a)}da\\
&=&T^{r-1}\int_{\mathfrak{a}^+_1\cap\{\alpha=0\}} e^{2T\rho(a)}da\ll T^{r-1}e^{\delta_\alpha T}=o(e^{(\delta-\varepsilon)T}).
\end{eqnarray*}
Since
$$
\mathfrak{a}^+_T-\mathfrak{a}^\eta_T\subset \bigcup_{\alpha\in\Sigma^+} \mathfrak{a}^+_T\cap\{\alpha<\eta\}.
$$
This proves part (a) of the lemma. Part (b) follows from (\ref{eq_cartan}).

To prove part (c), we note that
$$
\hbox{\rm Vol}(G_{T+\varepsilon})=\int_{\mathfrak{a}^+_{T+\varepsilon}}\xi(a)da=(T+\varepsilon)^r\int_{\mathfrak{a}^+_{1}}\xi((T+\varepsilon)a)da
$$
It is easy to check that there exist $b>0$ such that $\sinh(t+\varepsilon)\le e^\varepsilon\sinh (t)+b$ for every $\varepsilon\in (0,1)$ and $t\ge 0$.  Thus, for $a\in\mathfrak{a}^+_1$ and sufficiently small $\varepsilon>0$,
$$
\xi((T+\varepsilon)a)\le \prod_{\alpha\in\Sigma^+} (a_\varepsilon\sinh(\alpha(Ta))+b)^{m_\alpha}\le d_\varepsilon\xi(Ta)+C\sum_{i\in I} e^{\lambda_i(a)}
$$
where $d_\varepsilon\to 1$ as $\varepsilon\to 0^+$, $C>0$, and $\lambda_i$'s are characters such that $2\rho-\lambda_i<0$ in the interior of $\mathfrak{a}^+$.
Thus, it follows from (\ref{eq_lambda}) that
$$
\int_{\mathfrak{a}_T^+} \xi((T+\varepsilon)a)da\le  d_\varepsilon \int_{\mathfrak{a}_T^+} \xi(Ta)da +o(e^{(\delta-\varepsilon)T}).
$$
Using (\ref{eq_gt}) and (\ref{eq_low}), we deduce that
$$
\limsup_{T\to\infty}\frac{\hbox{\rm Vol}(G_{T+\varepsilon})}{\hbox{\rm Vol}(G_{T})}\le d_\varepsilon,
$$
and part (c) of the lemma follows. The last part of lemma can be proved similarly.
\end{proof}

\section{Acknowledgements}

The main ideas of this paper were developed during the workshop ``Ergodic properties of geometric group actions'' in Summer 2003.
The authors would like to express deep appreciation to the organizers of this workshop and to the Max Planck Institute of
Mathematics for its support. We also would like to thank to R.~Spatzier for raising the problem solved in this paper during the workshop
and to Y.~Guivarc'h for bringing the paper of Arnold and Krylov to out attention.

\end{document}